\documentclass[letterpaper, 10 pt, conference]{ieeeconf}  
\IEEEoverridecommandlockouts                              
\usepackage{centernot}
\usepackage[noadjust]{cite}
\usepackage{longtable}
\usepackage{amsmath}
\usepackage{graphicx, amssymb}
\usepackage[dvips]{epsfig}
\usepackage{color}
\usepackage[dvipsnames]{xcolor}
\usepackage{comment}
\usepackage{balance}
\usepackage{tikz}

\newtheorem{pro}{Proposition}
\newtheorem{deff}{Definition}

\newtheorem{ex}{Example}

\def\tt{\tilde \theta}

\def\L2{{\cal L}_2}

\newlength{\defbaselineskip}
\newcommand{\setlinespacing}[1]%
           {\setlength{\baselineskip}{#1 \defbaselineskip}}

\newcommand*{\Scale}[2][4]{\scalebox{#1}{$#2$}}%

\newcommand{\beq}{\begin{equation}}
\newcommand{\eeq}{\end{equation}}

\newcommand{\carre} {\hfill $\blacksquare$}

\makeatletter
\newcommand{\xleftrightarrow}[2][]{\ext@arrow 3359\leftrightarrowfill@{#1}{#2}}
\makeatother

\newcommand{\xdasharrow}[2][-->]{
\tikz[baseline=-\the\dimexpr\fontdimen22\textfont2\relax]{
\node[anchor=south,font=\scriptsize, inner ysep=1.5pt,outer xsep=8pt](x){#2};
\draw[shorten <=3.4pt,shorten >=3.4pt,dashed,#1](x.south west)--(x.south east);
}
}

\def\BibTeX{{\rm B\kern-.05em{\sc i\kern-.025em b}\kern-.08em
    T\kern-.1667em\lower.7ex\hbox{E}\kern-.125emX}}

\usepackage{longtable}
\usepackage{amsmath}
\usepackage{graphicx}
\usepackage{color}
\usepackage[dvips]{epsfig}
\usepackage{times}

\usepackage{graphicx, amssymb}
\usepackage{amsmath}

\smallskip

\newtheorem{theo}{Theorem}

\newtheorem{lem}{Lemma}

\allowdisplaybreaks

\title{\LARGE \bf
Strong Structural Controllability of Signed Networks
}

\author{Shima Sadat Mousavi$^{\dagger}$, Mohammad Haeri$^{\dagger}$, and Mehran Mesbahi$^{\ddagger}$
\thanks{$^{\dagger}$The authors are with the Department of Electrical Engineering, Sharif University of Technology, Tehran, Iran. Emails: shimasadat$\_$mousavi@ee.sharif.edu, haeri@sharif.ir.}%
\thanks{$^{\ddagger}$The author is with the Department of Aeronautics and Astronautics, University of Washington, WA 98195. Email: mesbahi@uw.edu.}%
}
\begin{document}

\maketitle
\thispagestyle{empty}
\pagestyle{empty}

\begin{abstract}
In this paper, we discuss the controllability of a family of linear time-invariant (LTI) networks defined on a signed graph. In this direction, we introduce the notion of positive and negative signed zero forcing sets for the controllability analysis of positive and negative eigenvalues of system matrices with the same sign pattern. A sufficient combinatorial condition that ensures the strong structural controllability of signed networks is then proposed. Moreover, an upper bound on the maximum multiplicity of positive and negative eigenvalues associated with a signed graph is provided.    
\end{abstract}


\section{Introduction}

Thanks to the ubiquity and wide recent applications of networks, there has been a surge of interest in studying networked dynamical systems and their control. 
One of the fundamental problems pertinent to the control of networks is their controllability~\cite{tanner2004controllability}. 
In most cases, the exact value of the entries of system matrices, that is, the connection weights of a network, is unknown or highly uncertain.  Accordingly, finding alternative means of system analysis based on topological features of the underlying graph is of importance; these features are also instrumental in network design problems~\cite{egerstedt2012interacting,shima}. There are different works in the literature, adopting diverse points of view towards controllability analysis of networks. In some works, controllability of a particular dynamics, e.g., Laplacian dynamics, has been considered \cite{tanner2004controllability,mousavi2018controllability, mousavi2018laplacian}, while in other works, instead of a specific dynamics, a family of dynamical networks all of which are defined on the same structure, has been studied. The second approach leads to structural controllability analysis for networked systems of interest in this work.

In the structural controllability framework, the network is viewed in terms of the zero-nonzero pattern  of  system matrices. 
 In this direction, 
strong structural controllability results provide conditions ensuring the controllability for \emph{all} systems with the same zero-nonzero pattern \cite{mayeda1979strong}. In the systems and control literature, different interpretations of strong structural controllability have been presented in terms of spanning cycle \cite{jarczyk2011strong}, constrained $t$-matchings \cite{chapman2013strong}, and zero forcing sets \cite{monshizadeh2014zero,trefois2015zero,mousavi2018structural, mousavi2017robust,mousavi2018null, mousavi2019strong}.

Signed networks have recently attracted a lot of attention in the systems community; the controllability of 
this class of networks with a particular Laplacian  dynamics has also been examined in a few works \cite{sun2017controllability,she2018controllability}. In fact, signed networks can be representative of a wide range of scenarios of practical interest, such as social networks and fault tolerant networks~ \cite{altafini2013consensus,wasserman1994social}. In a signed network, the graph admits both positive and negative edges that indicate respectively, cooperative or adversarial interactions among the nodes. 
As such, by considering a \emph{sign pattern} instead of a zero-nonzero pattern, not only can we capture the network structure, but also define a more restrictive family of networks that represents distinct qualitative features. 

 The notion of \emph{sign controllability}, which is strong structural controllability of networks with the same sign pattern, was first introduced in \cite{johnson1993sign} and examined for the special case of single-input systems with all nonzero entries as positive. These results were later extended to the multi-input case in \cite{tsatsomeros1998sign}, where signed networks are examined in the context of the so-called \emph{strict linear control systems} with some restrictive properties; for example, the diagonal entries of the system matrices should be nonzero and have the same sign. In \cite{tsatsomeros1998sign}, sufficient algebraic conditions for sign controllability of a network as well as necessary and sufficient conditions for the sign controllability of {strict linear control systems} have been presented; however, recognition of these algebraic conditions was proven to be NP-hard. More recently, in \cite{hartung2013characterization}, sign controllability of another family of networks has been analyzed, and an algebraic condition has been provided for systems whose sign pattern admits only real eigenvalues. However,  the verification of these conditions is also NP-hard.


The notion of zero forcing game,  played on  a graph to change the color of the nodes based on a coloring rule, was defined in \cite{work2008zero} for the minimum rank problem.  
 Later, other variants of the zero forcing sets were introduced in~\cite{brimkov2019computational}. For example, in \cite{goldberg2014zero}, in order to study the minimum rank problem for symmetric matrices with the same sign pattern (with an undirected graph), a signed zero forcing set was defined. In  this paper, we introduce the new notion of positive and negative signed zero forcing sets for a directed signed graph  that can be utilized in providing an upper bound on the maximum geometric multiplicity of positive and negative eigenvalues of matrices with the same sign pattern. 

As the main contribution of this work, using the notion of signed zero forcing sets, strong structural controllability conditions for the zero, positive, and negative eigenvalues of matrices with the same sign pattern are provided. Furthermore, we present a sufficient condition for the strong structural controllability of signed networks, whose sign patterns admits only real eigenvalues. However, there is no restriction on the sign of diagonal entries, and we allow the self-loops of the signed networks not to have any specified signs.  In \cite{eschenbach1991sign}, a complete characterization of such networks has been provided. For example, one can mention undirected networks with symmetric pattern matrices. 
  A few examples are used throughout the paper to better illustrate the results.

\section{Preliminaries}
We denote the set of real numbers by $\mathbb{R}$. 
For a vector $v$, $v_i$ is its $i$th entry; for a matrix $M$,  $M_{ij}$ is the entry in row $i$ and column $j$. A subvector $v_X$ is comprised of $v_i$, for $i\in X$, ordered lexiographically. We denote the transpose of the matrix $M$ by $M^T$. 
The $n\times n$ identity matrix is  denoted by $I_n$, and  its $j$th column  is designated by $e_j$. 
We designate  by $|S|$ the cardinality of the set $S$. The sign function $\mathrm{sign}(.):\mathbb{R}\rightarrow \{+,-,0\}$ returns the sign of a nonzero scalar, and we have $\mathrm{sign}(a)=0$ if and only if $a=0$. We also define the sign inversion function  as $\mathrm{inv}(+)=-$  and $\mathrm{inv}(-)=+$.

A \emph{zero-nonzero pattern} $P\in \{\times,0,?\}^{n\times n}$ is a  matrix whose off-diagonal entries can be zero or nonzero, respectively denoted by 0 or $\times$, and the diagonals are chosen from the set $\{\times,0, ?\}$. The zero-nonzero pattern of a matrix $A$ is a matrix $P$ such that for $i\neq j$, $P_{ij}= 0$ if and only if $A_{ij}= 0$. Note that if $P_{ii}=?$, $A_{ii}$ can be both zero or nonzero.

A \emph{sign pattern}  $P_s\in \{+,-,0,?\}^{n\times n}$ is a matrix whose off-diagonal entries are from the set $\{+,-,0\}$, and the diagonals belongs to $\{+,-,0,?\}$ . The sign pattern of a matrix $A$ is some  $P_s$ such that for $i\neq j$, $(P_{s})_{ij}=\mathrm{sign}(A_{ij})$; also, $(P_s)_{ii}=\mathrm{sign}(A_{ii})$ whenever $(P_s)_{ii}\neq ? $. If $(P_s)_{ii}=?$, $A_{ii}$ can be zero or a nonzero with a positive or negative sign. 

A \textit{graph} is denoted by $G=(V,E, P)$, where $V=\{1,\ldots,n\}$ is the vertex set and $E\subseteq V\times V$ is the edge set of the graph. We write $(i,j)\in E$ when there is an edge from the node $i$ to the node $j$. $P$ is a  zero-nonzero pattern such that $(i,j)\in E$ whenever $P_{ji}\neq 0$. Note that in our setup, a graph $G$ can contain self-loops as $(i,i)$  for some $i\in V(G)$; if we have $P_{ii}=?$ for some $1\leq i\leq n$, we assign a label $?$ to the self-loop $(i,i)\in E$,  implying that $(i,i)$ can appear or not appear in  $G$.  For $(i,j)\in E$, node $j$ (resp., node $i$) is an out-neighbor (resp., in-neighbor) of node $i$ (resp., node $j$). We denote by $N_{out}(i)$ 
the set of out-neighbors
 of node $i$. An undirected graph is a graph such that
$(i, j)\in E(G)$ if and only if $(j, i) \in E(G)$; in this case,
we write $\{i, j\} \in E(G)$, and node $j$ is referred to as the
neighbor of node $i$. The matrix $P$ is symmetric for an undirected graph.

A \textit{signed graph} $G_s$ is denoted by $G_s=(V,E, P_s)$, where 
 $P_s$ is an $n\times n$ sign pattern such that $(i,j)\in E$ whenever $(P_s)_{ji}\neq 0$. Then, to every edge $(i,j)\in E$, where $i\neq j$, we can assign  a label $+$ or $-$, that indicates  whether the weight of the connection between the nodes $i$ and $j$ is positive or negative. Moreover, a self-loop $(i,i)$, $i\in V$,  can be labeled with $+$, $-$, or $?$. This implies that the self-loop of node $i$ has a weight that can be positive, negative, or unspecified (including zero value). We can also define an undirected signed graph associated with a symmetric $P_s$. 

%

A \textit{looped graph}  is obtained from a graph by putting a self-loop on every node $v\in V$ that does not have a self-loop itself. Before precisely defining a looped graph, let us consider the indices $s$ and $r$ chosen respectively from the sets $\{\times, +,-,0,?\}$ and $\{+,-\}$. Now, we let the \emph{sign equations} as $?+s=?$,  $0+s=s$,  $r+\mathrm{inv}(r)=?$, $r+r=r$, and  $\times+\times=?$. One can verify theses equations by considering different scalars with the same pattern denoted by some $s\in \{\times, +,-,0,?\}$ and checking the pattern of the result. For example, by adding two positive (resp., negative) numbers, a positive (resp., negative) number is obtained, leading to $r+r=r$. On the other hand, by adding a positive and a negative number, the result may be positive, negative, or zero, implying that $r+\mathrm{inv}(r)=?$.  
Now, let $\mathcal{D}(\times)$, $\mathcal{D}(+)$, and $\mathcal{D}(-)$ be $n\times n$ diagonal pattern matrices whose diagonals are respectively, $\times$, $+$, and $-$. Using the sign equations,  for a given zero-nonzero pattern $P$ and a sign pattern $P_s$, let us define $P^{\times}=P+\mathcal{D}(\times)$, $P_s^+=P_s+ \mathcal{D}(+)$, and  $P_s^-=P_s+ \mathcal{D}(-)$.  Then, for a graph $G=(V,E,P)$, we define the \emph{looped graph} $G^{\times}=(V,E,P^{\times})$. Moreover, for a signed graph $G_s=(V,E,P_s)$, one has the \emph{positive looped graph} $G_s^+=(V,E,P_s^+)$ and the  \emph{negative looped graph} $G_s^-=(V,E,P_s^-)$.

\begin{ex} For the zero-nonzero pattern $P$ and the sign pattern $P_s$ defined as $$
P=\begin{bmatrix}
? & 0 & \times\\ 0 & \times & 0 \\ 0 & \times & 0 
\end{bmatrix}, \:\:\: P_s=\begin{bmatrix}
? & - & 0 & 0\\ 0 & - & + & 0\\ 0 & 0 & + & 0\\0 & 0 & + & 0
\end{bmatrix},$$
the graphs $G$ and $G_s$ in Figs. \ref{G} (a) and \ref{Gs} (a) can be represented. Also, the looped graph $G^{\times}$  is shown in Fig. \ref{G} (b), and the positive and the negative looped graphs $G_s^+$ and $G_s^-$ are respectively, depicted in Figs. \ref{Gs} (b) and (c). 
\begin{figure}[hbt]
\includegraphics[width=.35\textwidth]{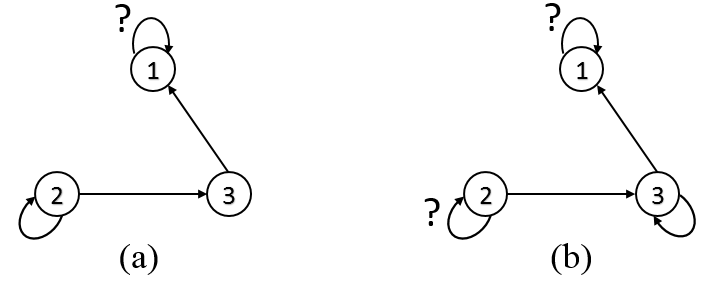}
\centering
\caption{a) Graph $G$, b) looped graph $G^{\times}$.}
\label{G} 
\end{figure}

\begin{figure}[hbt]
\includegraphics[width=.42\textwidth]{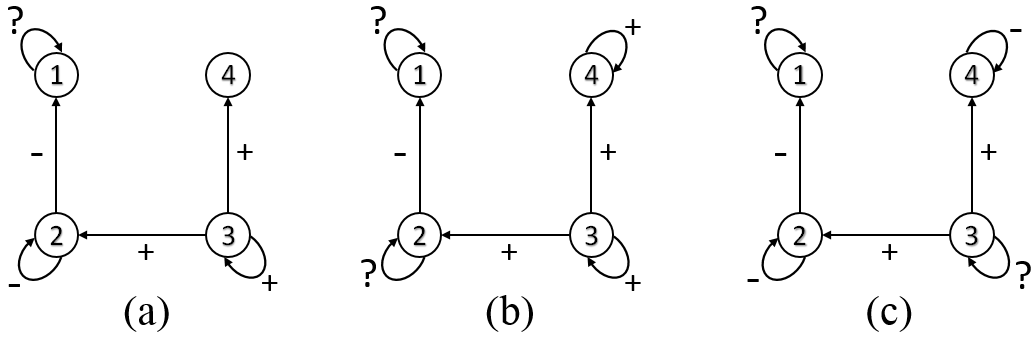}
\centering
\caption{a) Graph $G_s$, b) positive looped graph $G_s^{+}$, c) negative looped graph $G_s^{-}$.}
\label{Gs} 
\end{figure}
\end{ex}
 For a (undirected) graph $G=(V,E,P)$, the \emph{qualitative class}, denoted by $\mathcal{Q}(G)$, is defined as the set of all (symmetric) matrices in $\mathbb{R}^{n\times n}$ whose zero-nonzero pattern is $P$. Similarly, for a (undirected) signed graph $G_s=(V,E,P_s)$, the qualitative class $\mathcal{Q}_s(G_s)$, is the set of all (symmetric) matrices in $\mathbb{R}^{n\times n}$ whose sign pattern is $P_s$.
 
 We denote by $\Lambda(A)$ the set of eigenvalues of the matrix $A$. 
  For an eigenvalue $\lambda\in \Lambda(A)$, the dimension of the subspace $\mathcal{S}_{A}(\lambda)=\{\nu\in \mathbb{R}^n | \nu^TA=\lambda \nu^T \}$ is called the geometric multiplicity of  $\lambda$ and is  denoted by $\psi_{A}(\lambda)$. For some $\mathcal{M}\subseteq \Lambda(A)$, we also define the maximum geometric multiplicity of eigenvalues of $A$ belonging to $\mathcal{M}$ as $\Psi_{\mathcal{M}}(A)=\mathrm{max}\{\psi_{A}(\lambda)| \lambda\in \mathcal{M}\}$.

%
 
\subsection{Problem Formulation}
Given is an   LTI network with the following dynamics 
\begin{align}
\vspace{-.5em}
\dot{x}=Ax+Bu, 
\label{e1}
\end{align}  
where $x\in \mathbb{R}^n$ is the state vector of the nodes, and $u\in \mathbb{R}^m$ is the control input;
we refer to matrices $A\in \mathbb{R}^{n\times n}$ and $B\in \mathbb{R}^{n\times m}$ respectively, as the system and input matrices.
We let $A\in\mathcal{Q}_s(G_s)$ for some signed graph $G_s=(V,E,P_s)$; 
\color{black}{moreover}, $B$ is defined as
$B=[e_{j_1},\ldots,e_{j_m} ],$
\label{e2}
where nodes $j_k$, $k=1,\dots,m$, are called 
\textit{control nodes}, and  
 $V_C=\{j_1,\ldots,j_m\}$ is the set of control nodes. 

\color{black}{
  If with a suitable choice of the input, we can transfer the state of the nodes from any initial state to any final state within a finite time, then we say that the network with 
the pair $(A,B)$
 is controllable.}
\color{black}{As} controllability is preserved under a similarity transformation,
when the LTI system (\ref{e1}) is uncontrollable, 
there exists a nonsingular matrix $T\in \mathbb{R}^{n\times n}$ such that \color{black}{for some $q<n$,}
\vspace{-.5em}
\color{black}{\begin{equation}
\Scale[.87]{T^{-1}AT=\begin{bmatrix}\hat{A}_{11} & \hat{A}_{12}\\ 0 & \hat{A}_{22}\end{bmatrix}, \;
T^{-1}B=\begin{bmatrix}\hat{B}_1\\0\end{bmatrix}},
\label{decom}
\end{equation}}where $(\hat{A}_{11},\hat{B}_{1})$ is controllable, with $\hat{A}_{11}\in\mathbb{R}^{q\times q}$, $\hat{B}_{1}\in\mathbb{R}^{q\times m}$. 
 %
\color{black}{ When $\lambda\in \Lambda(A)$ and $\lambda\notin \Lambda(\hat{A}_{22})$, it is called a \textit{controllable eigenvalue} of the system (\ref{e1}). On the other hand, we define $\lambda$ as an \textit{uncontrollable eigenvalue} if  $\lambda\notin \Lambda(\hat{A}_{11})$. In this case, the input of the system cannot have any influence on $\lambda$. We can use \color{black}{the} Popov-Belevitch-Hautus (PBH) test for checking the controllability of eigenvalues. }     
 
\color{black}{\begin{pro}[\cite{sontag2013mathematical}]
The eigenvalue $\lambda$ of $A$ in a system with dynamics  (\ref{e1}) is controllable
if and only if for all nonzero $w$ for which $w^{T}A=\lambda w^T$, $w^{T}B\neq0$.
\end{pro}}
 \color{black}{An eigenvalue $\lambda$ is called  strongly structurally controllable  if it is a controllable eigenvalue for all $A\in \mathcal{Q}_s(G_s)$ for which $\lambda\in \Lambda(A)$.} \color{black}
Along the way, a \emph{signed} network with dynamics (\ref{e1}) which is defined on a signed graph $G_s=(V,E,P_s)$ is strongly structurally controllable if every $\lambda\in\Lambda(A)$ (for all $A\in\mathcal{Q}_s(G_s)$) is controllable. With a slight abuse of notation,  in this case, we say that $(G_s, V_C)$ is controllable.  Also, given a graph $G=(V,E,P)$, we say that $(G,V_C)$ is controllable  if  every $\lambda\in\Lambda(A)$ (for all $A\in\mathcal{Q}(G)$) is controllable.

\color{black}{Our focus is on the combinatorial characterizations
of strong structural controllability of positive, negative, and zero eigenvalues of a network, and then we provide a sufficient condition for strong structural controllability of signed networks whose sign patterns admits only real eigenvalues \cite{eschenbach1991sign}. 

\section{ Zero Forcing Games}

In this section, we first review the classical coloring rule and zero forcing sets for a graph $G$ \cite{barioli2009minimum}. Then the notions of \emph{signing  and coloring rule} and \emph{signed zero forcing sets} introduced in \cite{goldberg2014zero} are presented. Finally, the new notions of \emph{positive} and \emph{negative signed zero forcing sets} are discussed. The following
definitions can be utilized for undirected graphs by interpreting out-neighbors simply as neighbors.

\subsection{Classical Zero Forcing Sets}
Consider a graph $G=(V,E,P)$, and assume that some of its nodes are black, while the other nodes are white. The classical coloring rule is defined as follows.

\emph{\textbf{Classical coloring rule:}} Let $v\in V$ be either white with $P_{vv}\neq ?$ or black. If $v$ has only one white out-neighbor $u$, change the color of $u$ to black. 

Next, we define  classical and strong zero forcing sets. 

\begin{deff}
Assume that  $Z\subset V$ is the set of initially black nodes in   the graph $G=(V,E,P)$. The set $Z$ is a \emph{classical zero forcing set} of $G$ if by repetitively applying the classical coloring rule in $G$, all the nodes   become black.
\end{deff} 

\begin{deff} Consider the looped graph $G^{\times}=(V,E,P^{\times})$ associated with the graph $G=(V,E,P)$. We refer to a set $Z\subset V$ as a \emph{strong zero forcing set} of $G$ if by repetitively applying the classical coloring rule in $G^{\times}$, all of its nodes become black.
\end{deff} 

\begin{ex} Consider the graph $G$ in Fig. \ref{zfs}
(a) with the initial set of black nodes $Z=\{2,4,5\}$. First,  node 2 forces its only one white out-neighbor node 3 to be black, and then node 1 is forced by node 3 to be black. Since all nodes can be finally black through the successive application of the classical coloring rule, $Z$ is a classical zero forcing set of $G$. In addition, the looped graph $G^{\times}$ of the graph $G$ in Fig. \ref{nzfs} (a) is depicted in Fig. \ref{nzfs} (b). By performing the same chain of forces, one can see that $Z$ is a classical zero forcing set of $G^{\times}$, or equivalently, it is a strong zero forcing set of $G$.   
\begin{figure}[hbt]
\includegraphics[width=.42\textwidth]{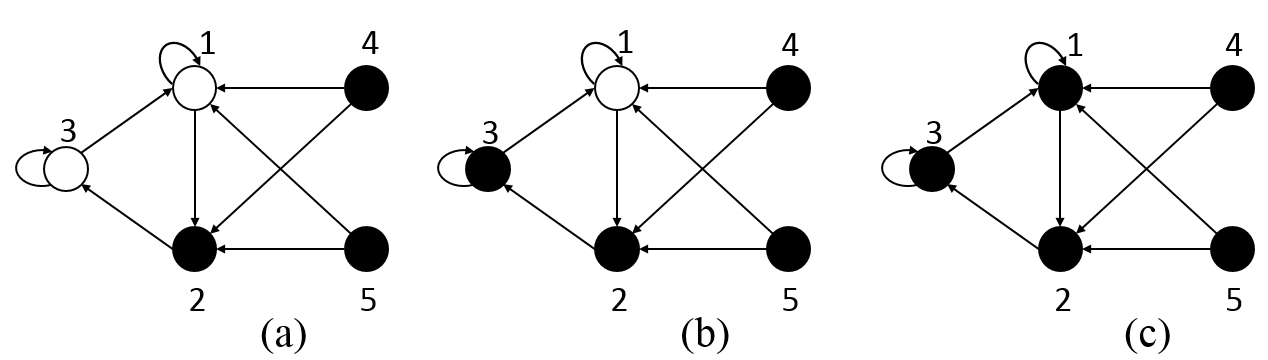}
\centering
\caption{An example for the classical coloring rule.}
\label{zfs} 
\end{figure}

\begin{figure}[hbt]
\includegraphics[width=.4\textwidth]{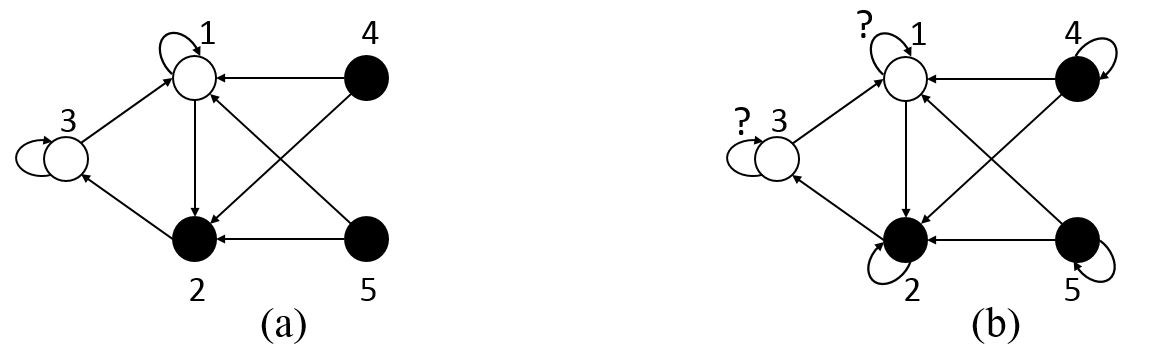}
\centering
\caption{a) Graph $G$, b) the associated looped graph $G^{\times}$.}
\label{nzfs} 
\end{figure}
\end{ex}

\subsection{Signed Zero Forcing Sets}

\emph{Signed zero forcing game} is indeed a signing and coloring game played on the nodes of a signed graph. In the first part of this game, we assume that some nodes of the signed graph $G_s=(V,E,P_s)$ are colored black, and others are white. Recall that a signed graph is a graph whose edges are labeled with the positive or  negative sign.  By doing this game, we aim to also assign the nodes of the graph a sign. For a node $u\in V$, let $m(u)$ denote its sign. If a node is assigned zero, its color is changed to black. Otherwise, if $u$ is white and is  marked with $+$ or $-$, we have $m(u)=+$ or $m(u)=-$. If a node is not marked, and its sign is undetermined, then we write $m(u)=*$. Thus, the goal of the game is to blacken the nodes and find the sign of white nodes when possible. 

Note that before starting the game, we only have some black nodes in the graph, and none of the nodes are marked with a sign. In this step, we can simply take one white node and mark it with $+$ and proceed based on the coloring rule. 
 
Before stating the game  rule, let us introduce some new notation. The letter $s$ is an index  taking values from $\{+,-\}$. If $s=+$, $\mathrm{inv}(s)=-$, and vice versa. For a node $v\in V$, let $W(v)=\{u\in N_{out}(v): u \:\: \mathrm{is} \:\: \mathrm{white}\}$. Then, $W(v)$ is the set of all white out-neighbors of the node $v$. Now, define $W_{+}(v)=\{u\in W(v): m(u)= (P_s)_{uv}\}$ and $W_{-}(v)=\{u\in W(v): m(u)=\mathrm{inv}( (P_s)_{uv})\}$. Accordingly, $W_{+}(v)$ (resp., $W_{-}(v)$) is the set of any white out-neighbor of the node $v$ which is marked, and its sign is the same as (resp., the opposite of)  the sign of the edge connecting $v$ to it. Also, let $W_*(v)=\{ u \in W(v): m(u)= *\}$. Then, $W_*(v)$ is the set of white out-neighbors of $v$ that has not yet been marked. 

Now, consider  a signed graph $G_s$ with all nodes colored either black or white, and some node of $G_s$ may be marked with $+$ or $-$. The rule of the game is stated as follows.

\emph{\textbf{Signing and coloring rule:}} Let $v\in V$ be either a black node or a white node with $(P_s)_{vv}\neq ?$ (then if $v$ is white, it has either no self-loops or a  self-loop labeled with $+$ or $-$). 
\begin{enumerate}
\item If $v$ has exactly one white out-neighbor $u$ (i.e. $W(v)=\{u\}$), then the color of $u$ is changed to black (note that $u$ and $v$ may be the same).
\item If   either $W_{+}(v)=W(v)$  or $W_{-}(v)=W(v)$, then all nodes in $W(v)$ become black. 
\item If all white out-neighbors of $v$ except one node $w$ are marked such that $W_{s}(v)\neq \emptyset$, $W_{\mathrm{inv}(s)}(v)=\emptyset$, and $W_*(v)=\{w\}$, then the unmarked node $w$ is marked with $P_{wv}.\mathrm{inv}(s)$.
\item If there is no white node in $G_s$ that is marked, and $u\in V$ is white, then  $u$ is marked with $+$.   
\end{enumerate}    
 
 Note that the first clause of the rule is the same as the classical coloring rule. 
 
In what follows, for a signed graph, the definitions of a signed zero forcing set, a positive signed zero forcing set,  and a negative signed zero forcing set are proposed.
 
 \begin{deff}
Let $Z\subset  V$ be a set of initially black nodes in the signed graph $G_s$. Apply the signing and coloring rules as many times as possible. The derived set of colored nodes of $Z$, denoted by $\mathcal{D}_c(Z)$, is defined as the final set of black nodes in $G_s$. Also, the derived set of marked nodes $\mathcal{D}_m(Z)$ is the set of any node $v$ with $m(v)=+\:\mathrm{or}\:-$ at the termination of the game. For an initial set of black nodes $Z$,  if $\mathcal{D}_c(Z)=V$, then $Z $ is called a \emph{signed zero forcing set} of $G_s$. 
 \end{deff}
 
  \begin{deff}
For a signed graph $G_s=(V,E,P_s)$,  consider the negative (resp., positive) looped graph $G^-_s=(V,E,P^-_s)$ (resp., $G^+_s=(V,E,P^+_s)$), and let $Z\subset  V$ be the set of initially black nodes. Now, perform the signing and coloring rule in $G^-_s$ (resp., $G^+_s$) as many times as possible.
The set of all nodes of $G_s$ that eventually become black in $G^-_s$ (resp., $G^+_s$) at the final stage of the game is called the positive (resp., negative) derived set of colored nodes of $Z$ and is denoted by $\mathcal{D}_c^+(Z)$ (resp., $\mathcal{D}_c^-(Z)$). Also, we denote the set of marked nodes at the termination of the game by $\mathcal{D}_m^+(Z)$ (resp., $\mathcal{D}_m^-(Z)$) and refer to it as the positive (resp., negative) derived  set of marked nodes of $Z$. Now, given an initial set of black nodes $Z$, if $\mathcal{D}_c^+(Z)=V$ (resp., $\mathcal{D}_c^-(Z)=V$), $Z$  is called a \emph{positive  signed zero forcing set} (resp., \emph{negative signed zero forcing set}) of $G_s$. The cardinality of  a positive (resp., negative) signed zero forcing set of $G_s$ is called the positive (resp., negative) signed zero forcing number, and is denoted by ${\mathcal{Z}_s^+}$ (resp., $\mathcal{Z}_s^-$). 

 \end{deff}
 
\begin{ex} 
 Consider the signed graph $G_s$ shown in Fig. \ref{szfs1} (a) where the nodes 4 and 5 are initially colored black. The different steps of applying the signing and coloring rule are shown in this figure. As shown in      
Fig. \ref{szfs1} (b), the 4th clause of the rule is applied, and node 1 is marked with $+$. Then, in Fig. \ref{szfs1} (c), we apply the 3rd clause of the rule for node 4 and mark its unmarked out-neighbor (node 2) with $-$. Next, the 2nd clause of the rule is performed for node 5, forcing nodes 1 and 2 to become black. Finally, the 1st clause of the rule is applied for node 2, and it forces its white out-neighbor to be black. Thus, since all nodes of the graph are finally black, the set $\{4,5\}$ is a signed zero forcing set of  $G_s$.  
\begin{figure}[hbt]
\includegraphics[width=.45\textwidth]{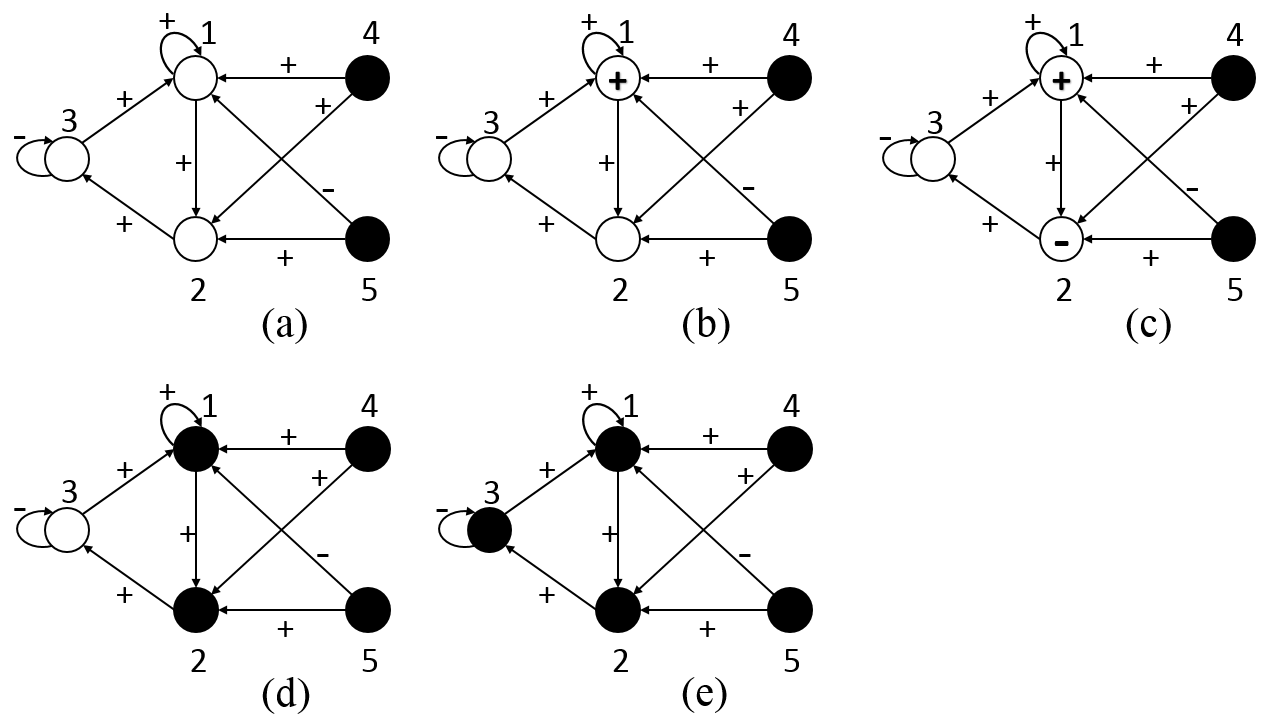}
\centering
\caption{An example of the signing and coloring rule. }
\label{szfs1} 
\end{figure}
Moreover, in Figs.  \ref{npszfs} (a) and (b), the negative looped graph $G_s^-$ and the positive looped graph  $G_s^+$ are respectively shown and through the application of a similar sequence of clauses of signing and coloring rule, we see that set $\{4,5\}$ is a signed zero forcing set of $G_s^-$ and $G_s^+$, and thus it is both a positive and a negative signed zero forcing  set of $G_s$.
%
\begin{figure}[hbt]
\includegraphics[width=.4 \textwidth]{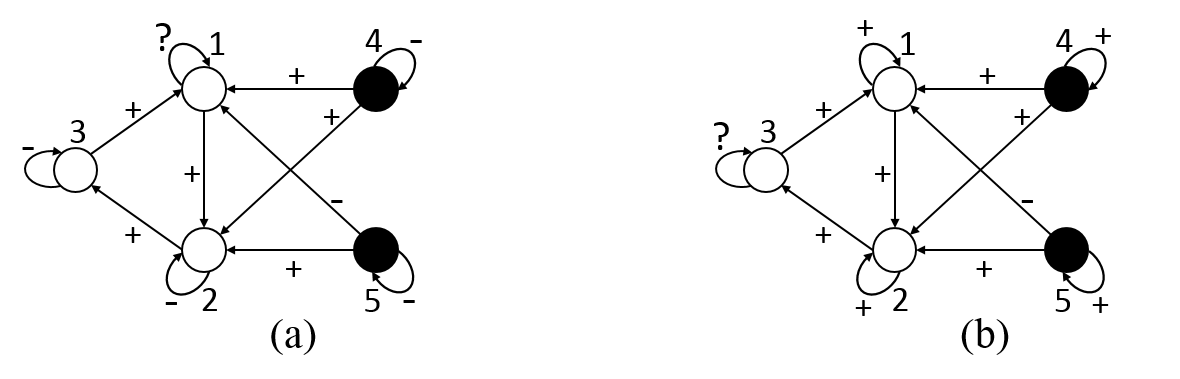}
\centering
\caption{a) Negative looped graph $G_s^-$, b) positive looped graph $G_s^+$  (associated with $G_s$ in Fig. \ref{szfs1}).}
\label{npszfs} 
\end{figure}
\end{ex}
\section{Strong Structural Controllability}
In this section, we derive combinatorial conditions, ensuring strong structural controllability of a signed network. 

The next theorem provides a necessary and sufficient condition for  controllability of $(G,V_C)$ in terms of classical zero forcing sets, where $G=(V,E,P)$. 

\begin{theo}[\cite{trefois2015zero}]
Given a network with dynamics (\ref{e1}), defined on the graph $G$, $(G,V_C)$ is controllable  if and only if $V_C$ is both a classical  and a strong zero forcing set of $G$. 
\end{theo}

\begin{ex} For  a network with dynamics (\ref{e1}) with the graph $G$ in Fig. \ref{zfs} (a), if $V_C=\{2,4,5\}$, $(G,V_C)$ is controllable, since we have shown that $V_C$ is both a classical and a strong zero forcing set of $G$. Note that set $V'_C=\{4,5\}$ cannot render the network strongly structurally controllable. This is due to the fact that although $V'_C$ is a classical zero forcing set of $G$, it is not a strong zero forcing set.  Indeed, for this network, the minimum number of control nodes ensuring the controllability of $(G, V_C)$ is 3.      
\end{ex}

Now, we consider a signed network and characterize the controllability  of its positive, negative, and zero eigenvalues in terms of the corresponding signed zero forcing sets. However, for the sake of brevity, we prove the results only for the positive eigenvalues; the proofs for the other cases are analogous.
  
\begin{lem}
 Let $G_s$ be a signed graph with a set of initially black nodes  $Z$ . 
 Let $A\in\mathcal{Q}_s(G_s)$, and $\nu\in \mathbb{R}^n$ be a left eigenvector of $A$ associated with $\lambda > 0$. If $\nu_i=0$ for all $i\in Z$, then $\nu_i=0$ for all $i\in \mathcal{D}_c^{+}(Z)$.
Moreover, if for some nodes $i\in \mathcal{D}_m^{+}(Z)$ with $\nu_i\neq 0$, one has $m(i)=\mathrm{sign}(\nu_i)$, then for any $k \in \mathcal{D}_m^{+}(Z)$ with $\nu_k\neq 0$,  $m(k)=\mathrm{sign}(\nu_k)$. 
 %
 \label{l2}
\end{lem}

\textit{Proof.}  Assume that the signed zero forcing game on the graph $G_s^+$ can be performed in $K$ step, at which only one of the clauses of the signing and coloring rule can be applied. In the first step, the nodes of $Z$ are colored black. Let $C_j$ and $M_j$ be respectively, the set of black nodes and marked nodes after the step $j$. Assume that $\nu_Z=0$. Note that for $j=1$, $C_j=Z$ and $M_j=\emptyset$. Also, $C_K=\mathcal{D}_c^{+}(Z)$ and $M_K=\mathcal{D}_m^{+}(Z)$. We claim that the theorem is not only true for $C_K$ and $M_K$, but also for any $C_j$ and $M_j$ that $1\leq j\leq K$. 
%
 The proof follows by a strong induction on $j$. It is clear that for $j=1$, the claim is true. 
 Now, assuming that the result holds for some $j\geq 1$, let us show its validity for $j+1$.  Consider the $i$th column of the matrix equation $\nu^TA=\lambda \nu^T$, that is, $\sum_{k\in N_{out}(i)}\nu_k A_{ki}= \lambda \nu_i$. This equation implies that 
\begin{equation}
\nu_i(A_{ii}-\lambda)+\sum_{k\in N_{out}(i)\setminus\{i\}}\nu_k A_{ki}= 0,
\label{eq3} 
\end{equation}          
where $\lambda>0$. Thus, if $\mathrm{sign}(A_{ii})\leq0$, we have $\mathrm{sign}(A_{ii}-\lambda)\leq0$, and otherwise $\mathrm{sign}(A_{ii}-\lambda)$ can be positive, negative, or zero. Then, it can be represented in the graph by a self-loop for node $i$ labeled with $-$ or $?$. Accordingly, we assume that for any node $v\in V$, $v\in N_{out}(v)$. Then, the matrix $A$ is replaced with a matrix $A^+\in P_s^+$.  Now, assume that the first clause is applied in step $j+1$. In other words, there is a node $v$ for which $(P_s)_{vv}\neq ?$ if it is white, and it has only one white out-neighbor $u$. Then, $C_{j+1}=C_{j}\cup \{u\}$. Also, equation (\ref{eq3}) leads to $\nu_u A^+_{uv}=0$, and as such, we have $\nu_u=0$, which shows the validity of the claim in this case. Regarding the 4th clause of the rule, note that when no white node is marked,  we can arbitrarily mark some white node $u$ with $+$; this follows from the fact that $\nu$ and $-\nu $ are both the eigenvectors, and any of which can be chosen in this case.  For the 2nd clause of the rule, equation (\ref{eq3}) simplifies to $\sum_{u\in W_+(v)}\nu_u A^+_{uv}=0$. Now, analogous to the proof of Theorem 3.2 in \cite{goldberg2014zero},
we claim that all summands on the right hand side of this equation have the same sign. Indeed, based on the definition, for all  $u\in W_+(v)$, we have $m(u) A^+_{uv}>0$. Moreover, according to the hypothesis of induction and without loss of generality, let us assume that for any $u\in W_+(v)$, $m(u)=\mathrm{sign}(\nu_u)$. Then, the claim immediately follows. Now, consider the case when the 3rd clause is applied. Equation (\ref{eq3}) in this case leads to $\sum_{u\in W_s(v)}\nu_u A^+_{uv} + \nu_w A^+_{wv}=0 $. Similar to the second case, we can prove that all summands of            $\sum_{u\in W_s(v)}\nu_u A^+_{uv}$ have the same sign. Without loss of generality, assume that for all $i\in W_s(v)$, $m(i)=\mathrm{sign}(\nu_i)$. Based on equation (\ref{eq3}), we should have $\mathrm{sign}(\nu_w)=P_{wv} \mathrm{inv}(s)$. Moreover, based on the 3rd clause of the rule, $m(w)=P_{wv} \mathrm{inv}(s)$, and hence we have $m(w)=\mathrm{sign}(\nu_w)$. Hence, the claim remains valid in this case, completing the proof. 
\carre

The next theorem is one of the main results of the paper, providing a sufficient condition for strong structural controllability of positive (resp., negative/zero) eigenvalues of a signed network. 

\begin{theo}
In an LTI network with  a signed graph $G_s$, every positive (resp., negative/zero) eigenvalue of all $A\in \mathcal{Q}_s(G_s)$ is controllable if $V_C$ is a positive  signed zero forcing set (resp., negative signed zero forcing set/signed zero forcing set) of $G_s$. 
\label{TH}
\end{theo}

\textit{Proof.} We only provide the proof for controllability of positive eigenvalues for brevity. Suppose $V_C$ is a  positive signed zero forcing set of $G_s$, but there is some $A\in\mathcal{P}_s(G_s)$ 
with an uncontrollable positive eigenvalue $\lambda$. 
 Then, there is a nonzero left eigenvector $\nu$ associated with $\lambda>0$  such that $\nu^TB=0$, or equivalently, $\nu_i=0$, for all $i\in V_C$. Since $\mathcal{D}^{+}_c(V_C)=V$, from Lemma \ref{l2}, we have  $\nu=0$; which is a contradiction. 
 %
 %
 %
 %
%
 \carre
 
 The next example is mentioned in \cite{tsatsomeros1998sign}.
 
 \begin{ex} Consider a signed network with dynamics (\ref{e1}) with the signed graph $G_s$ in Fig. \ref{szfs1} (a), and let $V_C=\{4,5\}$. Since $V_C$ is a signed, a positive signed,  and a negative signed zero forcing set of $G_s$, then based on Theorem \ref{TH}, zero, positive, and negative eigenvalues of the network  are strongly structurally controllable. 
 \end{ex} 
 
Theorem 2 leads to the next result, a sufficient condition for strong structural controllability of signed networks. 
 
 \begin{theo}
 Consider an LTI network with a signed graph $G_s$ whose sign pattern admits only real eigenvalues. Then, $(G_s,V_C)$ is controllable  if $V_C$ is a signed, a positive signed, and a negative signed zero forcing set of $G_s$.
 \label{th33} 
 \end{theo}

   \begin{ex} 
 For the undirected network shown in Fig. \ref{ex}(a), one can verify that $V_C$ is a signed, positive signed, and negative signed zero forcing set of the graph; as such based on Theorem \ref{th33}, we can deduce that the signed network is strongly structurally controllable. To demonstrate that $V_C$ is a signed zero forcing set, the steps of the signing and coloring rule are shown in Fig. \ref{ex}.
 
 \begin{figure}[hbt]
\includegraphics[width=.28\textwidth]{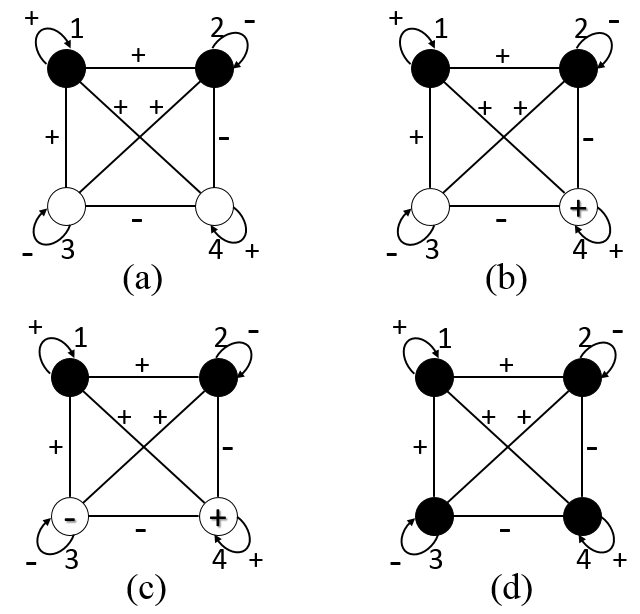}
\centering
\caption{An example of the signing and coloring rule. }
\label{ex} 
\end{figure}
\end{ex}

In \cite{goldberg2014zero}, an upper bound on the maximum nullity of matrices with a symmetric sign pattern has been obtained. In the same direction, by using the notions of positive and negative signed zero forcing sets, we can provide an upper bound on the maximum geometric multiplicity of positive and negative eigenvalues of matrices with the same sign pattern. Let $\Lambda_+(A)$ (resp., $\Lambda_-(A)$) denote the set of positive (resp., negative) eigenvalues of the matrix $A$. 
  
  \begin{pro}
Consider a signed graph $G_s=(V,E,P_s)$ with the positive (resp., negative) signed zero forcing number $\mathcal{Z}_s^+$ (resp., $\mathcal{Z}_s^-$). Then, for all $A\in \mathcal{Q}_s(G_s)$, we have $\Psi_{\Lambda_{+}(A)}(A)\leq~\mathcal{Z}_s^+$ (resp., $\Psi_{\Lambda_{-}(A)}(A)\leq~\mathcal{Z}_s^-$).    
  \label{p8}
  \end{pro}

  \textit{Proof:} 
Here we state the proof for positive eigenvalues only. 
 Assume  that for some $\lambda>0$, $\psi_A(\lambda)=k$. Then, similar to the proof of Proposition 2.2 of \cite{work2008zero}, we can say that for every subset of nodes  $X$ whose cardinality is $k-1$,  there exists a nonzero $\nu\in \mathcal{S}_{A}(\lambda)=\{\nu\in \mathbb{R}^n | \nu^TA=\lambda \nu^T \}$ such that $\nu_i=0$, for every $i\in X$. 
 Now, suppose that there exists a positive eigenvalue $\beta$ of some $A\in\mathcal{Q}_s(G)$ such that  $\psi_A(\beta)> \mathcal{Z}_s^+$.
 Let $Z$  be a positive signed zero forcing set of $G_s$ with the cardinality $\mathcal{Z}^+_s$. Then, there is a nonzero eigenvector $\nu\in \mathcal{S}_A(\beta)$ for which we have $\nu_i=0$, for all $i\in Z$. Additionally, we know from Lemma \ref{l2} that with the positive signed zero forcing set $Z$, if $\nu_Z=0$, then $\nu=0$, as $\mathcal{D}^+_c(Z)=V$. Thus, we reach a contradiction. \carre 

\balance 
\section{Conclusion}

Strong structural controllability of positive, negative and zero eigenvalues  of LTI systems defined on the same signed graphs has been examined in this work. We introduced the notions of positive and negative signed zero forcing sets; these notions can  be used to provide a set of control nodes for ensuring strong structural controllability of positive and negative eigenvalues of signed networks. Moreover, we have shown that a signed zero forcing set,  as a set of control nodes,  renders the zero eigenvalues strongly structurally controllable. Finally, an upper bound on the maximum multiplicity of the positive and negative eigenvalues of matrices with the same sign pattern has been presented. 

\bibliographystyle{IEEEtran}
\bibliography{library}

\begin{thebibliography}{10}
\providecommand{\url}[1]{#1}
\def\UrlFont{\rmfamily}
\providecommand{\newblock}{\relax}
\providecommand{\bibinfo}[2]{#2}
\providecommand\BIBentrySTDinterwordspacing{\spaceskip=0pt\relax}
\providecommand\BIBentryALTinterwordstretchfactor{4}
\providecommand\BIBentryALTinterwordspacing{\spaceskip=\fontdimen2\font plus
\BIBentryALTinterwordstretchfactor\fontdimen3\font minus
  \fontdimen4\font\relax}
\providecommand\BIBforeignlanguage[2]{{%
\expandafter\ifx\csname l@#1\endcsname\relax
\typeout{** WARNING: IEEEtran.bst: No hyphenation pattern has been}%
\typeout{** loaded for the language `#1'. Using the pattern for}%
\typeout{** the default language instead.}%
\else
\language=\csname l@#1\endcsname
\fi
#2}}

\bibitem{tanner2004controllability}
H.~G. Tanner, ``On the controllability of nearest neighbor interconnections,''
  in \emph{Proc. 43rd {IEEE} Conf. on Decision and Control}, vol.~3, 2004, pp.
  2467--2472.

\bibitem{egerstedt2012interacting}
M.~Egerstedt, S.~Martini, M.~Cao, K.~Camlibel, and A.~Bicchi, ``Interacting
  with networks: How does structure relate to controllability in single-leader,
  consensus networks?'' \emph{IEEE Control Syst. Mag.}, vol.~32, no.~4, pp.
  66--73, 2012.

\bibitem{shima}
S.~S. Mousavi and M.~Haeri, ``Controllability analysis of networks through
  their topologies,'' in \emph{Proc. 55th {IEEE} Conf. on Decision and
  Control}, 2016, pp. 4346--4351.

\bibitem{mousavi2018controllability}
S.~S. Mousavi, M.~Haeri, and M.~Mesbahi, ``Controllability analysis of
  threshold graphs and cographs,'' in \emph{2018 European Control Conf.}, 2018,
  pp. 1--6.

\bibitem{mousavi2018laplacian}
S.~S. Mousavi, M.~Haeri, and M.~Mesbahi, ``Laplacian dynamics on cographs:
  Controllability analysis through joins and unions,'' \emph{arXiv:1802.03599},
  2018.

\bibitem{mayeda1979strong}
H.~Mayeda and T.~Yamada, ``Strong structural controllability,'' \emph{SIAM J.
  Contr. and Optimiz.}, vol.~17, no.~1, pp. 123--138, 1979.

\bibitem{jarczyk2011strong}
J.~C. Jarczyk, F.~Svaricek, and B.~Alt, ``Strong structural controllability of
  linear systems revisited,'' in \emph{Proc. 50th {IEEE} Conf. on Decision and
  Control and Eur. Control Conf.}, Orlando, FL, 2011, pp. 1213--1218.

\bibitem{chapman2013strong}
A.~Chapman and M.~Mesbahi, ``On strong structural controllability of networked
  systems, a constrained matching approach,'' in \emph{Proc. American Control
  Conf.}, Washington, DC, 2013, pp. 6126--6131.

\bibitem{monshizadeh2014zero}
N.~Monshizadeh, S.~Zhang, and M.~K. Camlibel, ``Zero forcing sets and
  controllability of dynamical systems defined on graphs,'' \emph{IEEE Trans.
  Automat. Contr.}, vol.~59, no.~9, pp. 2562--2567, 2014.

\bibitem{trefois2015zero}
M.~Trefois and J.-C. Delvenne, ``Zero forcing number, constrained matchings and
  strong structural controllability,'' \emph{Linear Alg. and its Applic.}, vol.
  484, pp. 199--218, 2015.

\bibitem{mousavi2018structural}
S.~S. Mousavi, M.~Haeri, and M.~Mesbahi, ``On the structural and strong
  structural controllability of undirected networks,'' \emph{IEEE Trans.
  Automat. Contr.}, vol.~63, no.~7, pp. 2234--2241, 2018.

\bibitem{mousavi2017robust}
S.~S. Mousavi, M.~Haeri, and M.~Mesbahi, ``Robust strong structural
  controllability of networks with respect to edge additions and deletions,''
  in \emph{Proc. American Control Conf.}, 2017, pp. 5007--5012.

\bibitem{mousavi2018null}
S.~S. Mousavi, A.~Chapman, M.~Haeri, and M.~Mesbahi, ``Null space strong
  structural controllability via skew zero forcing sets,'' in \emph{2018
  European Control Conf.}, 2018, pp. 1845--1850.

\bibitem{mousavi2019strong}
S.~S. Mousavi, M.~Haeri, and M.~Mesbahi, ``Strong structural controllability
  under network perturbations,'' \emph{arXiv:1904.09960}, 2019.

\bibitem{sun2017controllability}
C.~Sun, G.~Hu, and L.~Xie, ``Controllability of multiagent networks with
  antagonistic interactions,'' \emph{IEEE Trans. Automat. Contr.}, vol.~62,
  no.~10, pp. 5457--5462, 2017.

\bibitem{she2018controllability}
B.~She, S.~Mehta, C.~Ton, and Z.~Kan, ``Controllability ensured leader group
  selection on signed multiagent networks,'' \emph{IEEE Trans. cybernetics},
  2018.

\bibitem{altafini2013consensus}
C.~Altafini, ``Consensus problems on networks with antagonistic interactions,''
  \emph{IEEE Trans. Automat. Contr.}, vol.~58, no.~4, pp. 935--946, 2013.

\bibitem{wasserman1994social}
S.~Wasserman and K.~Faust, \emph{Social Network Analysis: Methods and
  Applications}.\hskip 1em plus 0.5em minus 0.4em\relax Cambridge university
  press, 1994, vol.~8.

\bibitem{johnson1993sign}
C.~R. Johnson, V.~Mehrmann, and D.~D. Olesky, ``Sign controllability of a
  nonnegative matrix and a positive vector,'' \emph{SIAM J. Matrix Analysis and
  Applic.}, vol.~14, no.~2, pp. 398--407, 1993.

\bibitem{tsatsomeros1998sign}
M.~J. Tsatsomeros, ``Sign controllability: Sign patterns that require complete
  controllability,'' \emph{SIAM J. Matrix Analysis and Applic.}, vol.~19,
  no.~2, pp. 355--364, 1998.

\bibitem{hartung2013characterization}
C.~Hartung, G.~Reissig, and F.~Svaricek, ``Characterization of sign
  controllability for linear systems with real eigenvalues,'' in \emph{2013
  Australian Contr. Conf.}, 2013, pp. 450--455.

\bibitem{work2008zero}
{AIM Minimum Rank--Special Graphs Work Group}, ``Zero forcing sets and the
  minimum rank of graphs,'' \emph{Linear Alg. and its Applic.}, vol. 428,
  no.~7, pp. 1628--1648, 2008.

\bibitem{brimkov2019computational}
B.~Brimkov, C.~C. Fast, and I.~V. Hicks, ``Computational approaches for zero
  forcing and related problems,'' \emph{Europ. J. Operational Research}, vol.
  273, no.~3, pp. 889--903, 2019.

\bibitem{goldberg2014zero}
F.~Goldberg and A.~Berman, ``Zero forcing for sign patterns,'' \emph{Linear
  Alg. and its Applic.}, vol. 447, pp. 56--67, 2014.

\bibitem{eschenbach1991sign}
C.~A. Eschenbach and C.~R. Johnson, ``Sign patterns that require real, nonreal
  or pure imaginary eigenvalues,'' \emph{Linear and Multilinear Algebra},
  vol.~29, no. 3-4, pp. 299--311, 1991.

\bibitem{sontag2013mathematical}
E.~D. Sontag, \emph{Mathematical Control Theory: Deterministic Finite
  Dimensional Systems}.\hskip 1em plus 0.5em minus 0.4em\relax New York:
  Springer Verlag, 1998.

\bibitem{barioli2009minimum}
F.~Barioli, S.~M. Fallat, H.~T. Hall, D.~Hershkowitz, L.~Hogben, H.~Van~der
  Holst, and B.~Shader, ``On the minimum rank of not necessarily symmetric
  matrices: a preliminary study,'' \emph{Electron. J. Linear Algebra}, vol.~18,
  no.~1, pp. 126--145, 2009.

\end{thebibliography}

\end{document}